%% file: lub-Q-tensors.tex
\documentclass[12pt]{article}

\usepackage{amsmath}
\usepackage{amssymb}
\usepackage{rotate}
\usepackage{graphics}
\usepackage{color}
\usepackage{color} 
\usepackage{graphicx}
\usepackage[sort&compress,square,comma,numbers]{natbib}
\usepackage{seceqn}
\usepackage[a4paper]{geometry}
\newcommand{\wscale}{{\cal E}}
\usepackage{soul}
\usepackage[colorlinks=true]{hyperref}  
\input{macros-cell}

\date{}

\begin{document}  


\title{Thin film models for an active gel\\[1cm]}

\author{
Georgy Kitavtsev\thanks{E-mail: georgy.kitavtsev@bristol.ac.uk} $^{1}$
\and \hspace*{-1.5cm}
Andreas M\"unch\thanks{E-mail: muench@maths.ox.ac.uk} $^{2}$
\and \hspace*{-1.5cm}
Barbara Wagner\thanks{E-mail: wagnerb@wias-berlin.de} $^{3}$
\\
\small{$^1$School of Mathematics, University of Bristol, University Walk, Bristol BS8 1TW, UK}
\\
\small{$^2$Mathematical Institute, University of Oxford, Andrew Wiles Building,  Oxford OX2 6GG, UK}\\
\small{$^3$Weierstra{\ss} Institute, Mohrenstrasse 39, 10117 Berlin, Germany} \\
} 

\maketitle

\vspace{1cm}

\allowdisplaybreaks

\begin{abstract}
In this study we present a free-boundary problem for an active liquid crystal based on the Beris-Edwards theory that uses a tensorial order parameter and includes active contributions to the stress tensor to analyse the rich defect structure observed in applications such as the Adenosinetriphosphate (ATP) driven motion of a thin film of an actin filament network. The small aspect ratio of the film geometry allows for an asymptotic approximation of the free-boundary problem in the limit of weak elasticity of the network and strong active terms. The new thin film model captures the defect dynamcs in the bulk as well as wall defects and thus presents a significant extension of previous models based on the Lesli-Erickson-Parodi theory. Analytic expression are derived that reveal the interplay of anchoring conditions, film thickness and active terms and their control of transitions of flow structure. 

\end{abstract}

\section{Introduction}

Since the works by Simha and Ramaswamy\cite{Simha2002} and Kruse et al. \cite{KJJPS04} active liquid crystals have been used extensively as a hydrodynamic theory to describe the ordered motion of large numbers of self-propelled particles, such as bacterial suspensions, fibroblast monolayers, or the ATP driven actin network that underlies the movement of the lamelopodium of a crawling cell. The different levels of description, from the microcopic to the continuum hydrdynamic theory of this rapidly expanding research field has been reviewed in Marchetti et al. \cite{Marchetti2013}.

In many studies active matter extensions are based on the Leslie-Ericksen-Parodi (LEP) theory \cite{leslie79,Ericksen1962} such as in \cite{KJJPS05,Juelicher2007}, where active polar gels were derived from thermodynamic principles. As in passive liquid crystals, defects are a common phenomenon and their dynamics is strongly influenced by the fact that the system is out-of-equilibrium due to the energy source from the active terms. 
Observations in in-vitro experiments \cite{Nedelec1997} show that they may directly depend on strength of the activity, where it was demonstrated that the observed defects tend to disappear again for sufficiently high levels of activity \cite{Backouche2006}. 

Based on the LEP theory, point defects such as asters, vortices and spirals were described \cite{KJJPS05,KJJPS04}. Furthermore, phase diagrams of unbounded two-dimensional states \cite{VJP06} as well as flow transitions in confined films \cite{VJP05} were investigated. In particular, it was found that spontaneous flow arises in a confined active polar gel (with no-slip or free-slip conditions at the domain walls) above a critical layer thickness. This transition was also described within a thin-film model with a free, capillary surface \cite{Sankararaman2009}.  

However, there are some inherent deficiencies to desribe the complete defect structure of passive liquid crystals based on the Leslie-Ericksen-Parodi theory, which is connected to the discontinuity of the director field and the infinite associated local elastic energy at the defect points. This problem becomes even more critical for the description of wall and line defects along which the elastic energy in the Leslie-Ericksen-Parodi theory is essentially discontinuous and, in particular, standard energy renormalization techniques can not be applied.  Moreover, when modelling the evolution of thin nematic films with moving contact lines using LEP theory, related problems occur due to singularity of the director field at the contact line  \cite{AC01,cummings04,Cummings2011,Lin2013,Lin2013a,Manyuhina2010,Manyuhina2013,Lam2015}.
Therefore, more general approaches such as the Beris-Edwards theory \cite{Beris1994,Za12} of liquid crystal hydrodynamics, that use a tensorial order parameter, the so-called Q-tensor, instead of a director field, have been devised. Extensions of this theory by active terms go back to Marenduzzo et al. \cite{Marenduzzo2007a,Marenduzzo2007} and have been extended in two- and three dimensions to various problems involving different geometries, such as spherical shells \cite{Giomi2012,Giomi2014, Thampi2014b,Marbach2015,Thampi2016,Bonelli2016}. But even for the passive Beris-Edwards theory, the conditions at boundaries and in particular, stress and anchoring conditions at free interfaces are less well studied within this model. Important contributions to these issues can be found in \cite{rey00,Rey1999,Lin2013}. In particular, it was conjectured~\cite{Re07} that a $Q$-tensor based approach might facilitate the resolution of nematic point defects in the vicinity of moving film contact lines~\cite{CM09}. 

The derivation of the corresponding thin-film model is the goal of this article.
We begin by formulating the active Beris-Edwards model (section~\ref{sec:form}) including all the boundary conditions for a two-dimensional cross section of a thin film. We emphasize that the two-dimensional Beris-Edwards model resembles basic features of its full three-dimensional version, but also finds independent interesting applications for modelling biological films on curved surfaces~\cite{Sanchez2012, Giomi2014, RVK12, PZ12}. In this case, we are able to represent of the $Q$-tensor variable through a scalar order parameter $q$ and the director field $n$ and reduce the active Beris-Edwards model to the corresponding active Ericksen model~\cite{Za12} describing the evolution of $q$ and $n$. Making use of the scale separation of the thin-film geometry, a leading order approximation is derived (section~\ref{sec:lub-ericksen}) in the limit of weak elasticity and strong active terms to arrive at a new thin-film model, both for the passive and active cases.
We also show that our model formally reduces to the one based on Leslie-Ericksen-Parodi theory, when the scalar order parameter $q$ is homogeneous, and coincides with one of~\cite{Lin2013} in the passive case.

Finally, we derive explicit solutions for special cases of  flat constant films  and small angle mismatch between the anchoring conditions. They show that in the passive case a solution with non-homogeneous nematic field exists when certain relations between film thickness and nematic boundary conditions are satisfied. In the active case, this solution also demonstrates nonzero flow and can be spontaneously initiated from the homogeneous one, for example by increasing the film thickness, similar to the effect observed in~\cite{VJP05,VJP06}.

A discussion of further extensions and applications concludes the paper (section~\ref{sec:conc}). In Appendix A we present the rescaled Ericksen model under the thin-film approximation. In Appendix B we derive the polar thin-film model based on the Leslie-Ericksen-Parodi system with active terms.

\section{Formulations of active liquid crystals}\label{sec:form}
\subsection{Beris-Edwards model for an active gel}
The model in \cite{KJJPS05} and the simplified version in \cite{VJP05} can be
viewed as based on the framework of liquid crystal theory augmented by sources
of energy due to ATP hydrolysis that drives the system and makes the bulk of
the cell an {\em active} (polar) gel.  The bulk liquid i.e.\ the gel is
characterized mainly by the velocity and the director field, which describes
the averaged orientation of the actin filaments at a given point in space and
time.  The driving force is provided in their models via the chemical potential
difference of ATP and its hydrolysis products. This hydrolysis of ATP fuels the
molecular motors (and is also used for the polymerization and depolymerization
of the actin filaments).  Instead of treating the chemical potential
difference as a local quantity or a fixed constant as is done for the Leslie-Ericksen-Prodi formulation \cite{leslie79,Leslie1968,Ericksen1962}, we include the corresponding terms into the Beris-Edwards theory that uses Q-tensors and is popular in the liquid crystal literature
\cite{Majumdar2010,DOY01,Za12,IXZ15,ZZRP16} nd the recent overview \cite{Majumdar2017,JBall2017} as a more general alternative theory for liquid crystals. 
In a subsequent step, we express the
Q-tensor in terms of the director field and an additional scalar order
parameter to obtain the Ericksen model, for which we then derive the thin film
model in section~\ref{sec:lub-ericksen}.  An active gel model in terms of Q-tensors and its
subsequent reformulation is also given in \cite{Marenduzzo2007}, but we also need to include
appropriate conditions at the free interface, which we base on \cite{Rey99,rey00}.  

We only consider two-dimensional models here and introduce
a spatial domain $\Omega$ with coordinates $(x_1, x_3)$, while $t$ represents
time,  (see Fig. 1 for a schematic sketch of the geometry and variables involved).
\begin{figure}
\centering
\includegraphics*[width=0.58\textwidth]{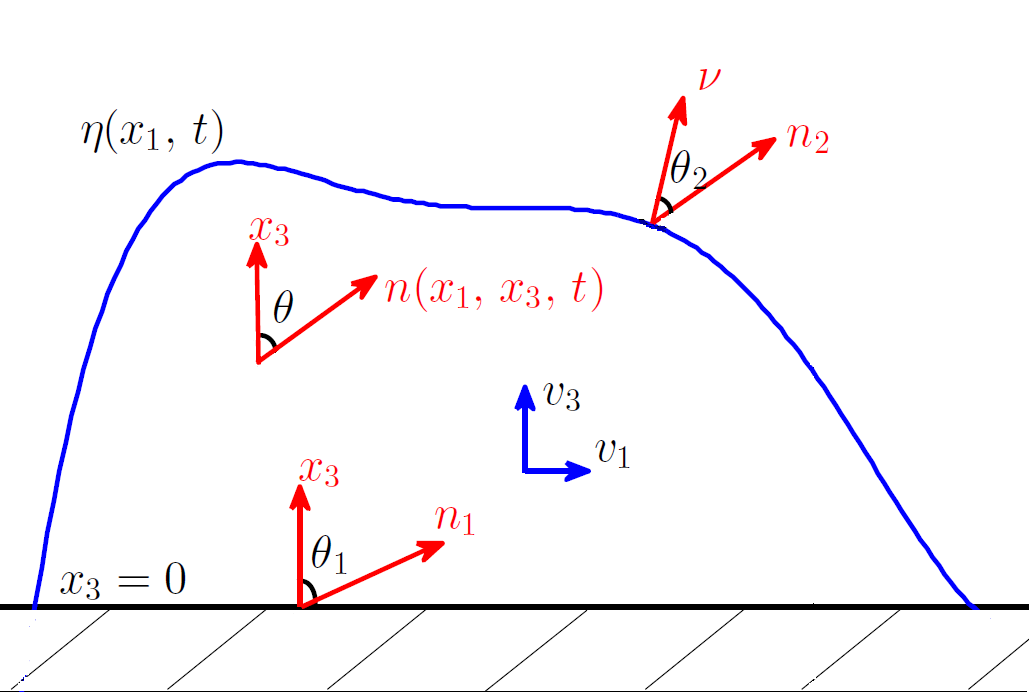}
\caption{ Sketch of the geometry of a thin film together with variables involved in the $Q$-tensor system \rf{cont1}-\rf{mol1}.} 
\label{Fig1}
\end{figure}
The Beris-Edwards model is associated with the standard Landau-de Gennes energy 
in the form \cite{KRSZ16}
\be
\mathcal{F}_{LG}[Q;\Omega]=\int_\Omega (f_e(Q)+f_b(Q))\,dx,
\lb{LGEn}
\ee
where $Q\in H^1(\Omega,\,\mathcal{L}_0)$ takes on values in the space of 
the symmetric and traceless matrices, or $Q$-tensors,
\bes
\mathcal{L}_0:=\{Q\in\mR^{2\times 2},\, Q=Q^T,\ \tr(Q)=0\}.
\ees
In \eqref{LGEn},
the bulk contribution is given by
\be
f_b(Q)=-\frac{a^2}{2}\tr{(Q^2)}+\frac{c^2}{4}\left(\tr{(Q^2)}\right)^2
\lb{bEn}
\ee
with $c>0$ and 
\be
f_e(Q)=\frac{L_1}{2}Q_{ij,k}Q_{ij,k}
\lb{ElEn}
\ee
is the elastic contribution, with an elastic constant $L_1>0$. 
(We deliberately avoid further complications by consider a model with 
only one elastic constant.)
Here, and elsewhere, we use the usual convention that duplicate indices are summed over
and indices with commas indicate spatial derivatives, e.g.\ $Q_{ij,k}$ is used for
the derivative of $Q_{ij}$ with respect to $x_k$. 

In the most general form, the Beris-Edwards model can be written as
(see e.g.~\cite{DOY01,Za12,IXZ15,ZZRP16} and references therein):
\begin{align}
0 &=\pai v_i, \label{cont1}\\
0 &=-\pai p +\mu\pj^2 v_i+\pj(\tau_{ij}+\sigma_{ij}-{\zeta\Delta\chi Q_{ij}})
		\label{Qmom1} \\
Q_t+(v\cdot\nabla)Q &= \Gamma H+S(\nabla v,Q)+{\lambda_1\Delta\chi\,Q}, 
\label{mol1}		
\end{align}  
where 
$v_i$ and $p$ are the velocity components and the pressure, and $\mu$ the 
isotropic viscosity.
The term
\be
S(\nabla v,Q)=(\xi e+\omega)(Q+\id/2)+(Q+\id/2)(\xi e-\omega)-2\xi(Q+\id/2)\tr(Q\nabla u),
\ee
with
\be
e_{ij}=\frac12\left(\pj v_i + \pai v_j\right),\quad\quad 
\omega_{ij}=\frac12\left(\pj v_i - \pai v_j\right),\quad I_{ij}=\delta_{ij},
\label{strain}
\ee
describes how the flow gradient rotates and stretches the order-parameter.  
The scalar parameter $\xi$ appearing both in equations \rf{Qmom1}  and \rf{mol1} depends on the molecular details of a given liquid crystal and measures the ratio between the tumbling and the aligning effect that a shear flow exert over the liquid crystal directors.
The active terms are associated with the activity parameters $\zeta$ and $\lambda_1$ 
and have been introduced in \rf{Qmom1}  and \rf{mol1} as in
\cite{ZZRP16,Marenduzzo2007}.
The
molecular field $H$ in \eqref{mol1} is the first variation of the Landau-de
Gennes energy \rf{LGEn} with respect to $Q$,
\begin{align}
H_{ij}=\frac{\delta {\cal F}_{LG}}{\delta Q_{ij}}
&=
a^2Q_{ij}-c^2Q_{ij}\tr(Q^2)+\lambda(x)\delta_{ij}
+L_1\partial_k^2 Q_{ij}.
\lb{H}
\end{align}
The Lagrange multiplier arises from the constraint $\mathrm{tr}\,Q=0$. We note that this 
constraint is equivalent to the normalisation condition of the director field,
as can be seen by taking the trace of the representation \rf{Qrep} for $Q$.
However, taking the trace of equation \rf{mol1}, gives, after some algebra, that
$
\lambda(x)=-2(\xi+1)Q_{il}\omega_{li}=0,
$
where the last equality follows from $Q_{il}=Q_{li}$. 
We will therefore drop the $\lambda(x)\id$ term from \eqref{H}.
The  symmetric and antisymmetric parts of the stress tensor $\sigma_{ij}$ and
$\tau_{ij}$ that appear in \eqref{Qmom1} are due to the director-flow
interaction and have the form
\begin{align}
\tau_{ij}&=-\xi(Q_{ik}+{\delta_{ik}}/{2})H_{kj}-\xi H_{ik}(Q_{kj}+{\delta_{kj}}/{2})
\notag\\
&\qquad +2\xi(Q_{ij}
+{\delta_{ij}}/{2})H_{km}Q_{km}-L_1\pj Q_{km}\pai Q_{km}
\lb{tau}
\end{align}
and
\be
\sigma_{ij}=Q_{ik}H_{kj}-H_{ik}Q_{kj}.
\lb{sigma}
\ee
For future reference, we also introduce the total stress tensor $T$, which includes
all contributions, including those from the active term, that is
\begin{equation}\label{Q:totalT}
T_{ij}=-p\delta_{ij}+2\mu e_{ij}+\tau_{ij}+\sigma_{ij}-\zeta\Delta\chi Q_{ij}.
\end{equation}

\paragraph{Boundary conditions at the substrate.}
We assume that the substrate is impermeable and that the no-slip condition
holds for the liquid, hence both components of the liquid vanish at $x_3=0$,
\begin{subequations}\label{Q:bcsubs}
\begin{equation}
v=0.
\end{equation}
We also impose strong anchoring, so that at $x_3=0$, we have
\begin{equation}
Q=Q_1=q_1(n_1\otimes n_1-\frac{1}{2}\id),
\end{equation}
\end{subequations}
with a given constant $q_1\in\mR$ and $n_1=[\sin(\theta_1),\,\cos(\theta_1)]\in\mRtwo$ (see also \cite{mottram2014}).

\paragraph{Boundary conditions at the free interface.}
We use the isotropic surface energy from~\cite{Rey99}
(retaining only the first constant term),
\be
F_s(Q,\nu)=g_0,
\lb{Sen}
\ee
which leads to the surface stress  (with $\id_s \equiv I-\nu\otimes\nu$)
\begin{equation}\label{Q:Ts1}
T^s=F_s \id_s
\end{equation}
that appears
in the right hand side of stress condition at the interface $x_3=\eta(x_1,t)$ 
\be\label{Q:Ts2}
\nu_i T_{ij}=(\delta_{ik}-\nu_i\nu_k)\pk T_{ij}^s.\\
\ee
In addition, we have the kinematic condition  
\begin{equation}
\eta_t=v_3-v_1\eta_{,1}
\end{equation}
at $x_3=\eta(x_1,t)$
and we impose the conical anchoring condition on $Q$, see also \cite{Teixeira1993,Ryschenko1976,JBall2017},
\begin{equation}
Q=q_2\left(R(\theta_2)\nu\otimes R(\theta_2)\nu-\frac{1}{2}\id\right),
\lb{con_bc}
\end{equation}
with a given constant $q_2\in\mR$ and $\theta_2\in [0,\pi)$, where
\bes
R(\theta)=\left[\begin{array}{ll}
\cos\theta & \sin\theta,\\
-\sin\theta & \cos\theta
\end{array}\right]
\ees
is the rotation by angle $\theta_2$.
We note, however, that in the thin film limit, the normal to the free boundary
is, to leading order, equal to the canonical unit vector $e_3$, hence this
boundary condition reduces to a strong anchoring condition with a fixed angle
$\theta_2$ with respect to the $x_3$ coordinate direction.

\subsection{Reduction to an active Ericksen model}
The reduction of the model \rf{cont1}--\rf{mol1} proceeds as follows:
By definition the two eigenvalues of $Q$  are $\pm q/2$
for some scalar order parameter $q\in\mR$. 
Moreover, one can show that for each $Q\in\mathcal{L}_0$ there exists a unit vector $n\in S^1$ (called {director}) such that representation
\be
Q=q[n\otimes n-\id/2]
\lb{Qrep}
\ee
holds. From this it also follows that each two-dimensional Q-tensor on a plane
is completely characterized by two degrees of freedom: the order parameter $q$
and the director  $n$. The representation \eqref{Qrep} does not distinguish
between $+n$ and $-n$. For definiteness, we fix the sign at the free interface, and hence
by continuity everywhere in the film, by requiring that $n$ points out of the
liquid and the director field is continuous everywhere in the film
bulk. In section 5 we will describe situations when the reduction presented in
this section can be extended without changes to the case of singular
director fields $n$ having defects.

We note also that under representation \rf{Qrep} the bulk energy \rf{bEn}
reduces to
\be
f_b(Q)=f_b(q)=-\frac{a^2q^2}{8}+\frac{c^2q^4}{64},
\ee
which attains its global minima at $q_{min}=\pm 2a/c$.

Substituting \rf{Qrep} into \rf{H} (and taking into account that $\lambda=0$) one obtains
\beas
H_{ij}&=&\left(a^2q-\frac{c^2q^3}{2}\right)\left[n_in_j-\frac{\delta_{ij}}{2}\right]
+L_1q_{,kk}\left[n_in_j-\frac{\delta_{ij}}{2}\right]\\
&&+2L_1q_{,k}n_{i,k}n_j+2L_1q_{,k}n_in_{j,k}+L_1q[2n_{i,k}n_{j,k}+n_{i,kk}n_j+n_in_{j,kk}]
\eeas

On the other hand, expressing of $H$ from \rf{mol1} gives
\be\label{E:GHij}
\Gamma
H_{ij}=q(n_jN_i+n_iN_j)+(q_t+v_kq_{,k})\left[n_in_j-\frac{\delta_{ij}}{2}\right]-S(\nabla v,Q)
-\lambda_1\Delta\chi q(n_i n_j-\delta_{ij}/2)
\ee
where we denote the rate of change of the director with respect to the background fluid
\be
N_i=\dot n_i -\omega_{ij} n_j,\quad\quad \dot n_i = \pt n_i + v_j\pj n_i,
\label{Ni}
\ee
and $\dot n_i$ denotes the material derivative.
 
Calculating the variational quantity $\Gamma(H_{ij}n_j+n_iH_{ij})$ for
both of the last two representations for $H$ and subsequently equating them
one obtains the following equation:
\begin{align}
&L_1\Gamma[2qn_{i,kk} {-2q|n_{i,k}|^2n_i}+q_{,kk}n_i+4q_{,k}n_{i,k}]+\Gamma(a^2q-\frac{c^2q^3}{2})n_i
\notag\\
=&
2qN_i+(q_t+v_kq_{,k})n_i-\frac{2}{3}(q+2)\xi e_{ji}n_j -\lambda_1\Delta\chi qn_i.
\lb{n_equ}
\end{align}
Multiplying the last equation by $n_i$ and using relations $N_in_i=n_i^2-1=0$
one obtains an Allen-Cahn type equation for the scalar order parameter
\be
q_t+v_kq_{,k}-\frac{2}{3}(q+2)\xi e_{ji}n_jn_i=L_1\Gamma q_{,kk} 
-4qL_1\Gamma|n_{j,k}|^2+\Gamma\left(a^2q-\frac{c^2q^3}{2}\right)
{+\lambda_1\Delta\chi q},
\lb{q_equ}
\ee
Using \rf{q_equ} one can simplify \rf{n_equ} to a parabolic equation for the
director field $n(x)$:
\be
L_1\Gamma[2qn_{i,kk}+4q_{,k}n_{i,k}]=2qN_i
-2qL_1\Gamma|n_{j,k}|^2n_i-\frac{2}{3}(q+2)\xi[e_{ji}n_j-e_{lk}n_ln_kn_i]
.
\lb{n_equ1}
\ee

Finally, the expressions \rf{sigma} and \rf{tau} for the symmetric and
antisymmetric stresses become
\bes
\Gamma\sigma_{ij}=q^2(n_iN_j-N_in_j)-\frac{\xi q(q+2)}{3}(n_in_ke_{kj}-e_{ik}n_kn_j)
\ees
and
\begin{align*}
\Gamma\tau_{ij}=&
-\frac{q\xi}{3}(q+2)(n_jN_i+n_iN_j)
+\frac{q\xi^2}{3}(4-q)(e_{ik}n_kn_j+n_in_ke_{kj})
\\ &\quad
+\frac{2\xi^2}{3}(q-1)^2e_{ij}
-\frac{8q^2\xi^2}{3}(\frac{3}{4}+q-q^2)\xi n_in_j
e_{ik}n_ln_k
\\ &\quad
+\frac{\xi q}{2}n_in_j(q_t+v_kq_{,k})
-\Gamma\,L_1\left(\frac{3}{4}q_{,i}q_{,j}+2q^2n_{k,i}n_{k,j}\right)
\\ &\quad
+
\xi\lambda_1\Delta \chi (1-q^2)q(n_i n_j-\delta_{ij}/2), 
\end{align*}
where the last term appears upon inserting the expression  \eqref{E:GHij} for
$H_{ij}$ into \eqref{tau}. Finally, we also have the explicit appearance of
the active stress in \eqref{Qmom1}, 
\begin{align*}
-\zeta\Delta\chi Q_{ij} &= -\zeta\Delta\chi q (n_i n_j-\delta_{ij}/2),
\end{align*}
so that the total stress tensor \eqref{Q:totalT} becomes
\begin{align}
T_{ij}&=-p\delta_{ij}+ T^E_{ij}  +   \tilde T_{ij} 
,
\label{E:totalT}
\intertext{with}
T^E_{ij} &= -L_1\left(\frac{3}{4}q_{,i}q_{,j}+2q^2n_{k,i}n_{k,j}\right),
\\
\tilde T_{ij} &= \alpha_1 n_kn_p e_{kp}n_in_j + \alpha_2 N_in_j + \alpha_3N_jn_i 
\notag\\&\quad
+\alpha_4 e_{ij} + \alpha_5 e_{ik}n_kn_j + \alpha_6 e_{jk}n_kn_i
+\frac{\xi q}{2\Gamma}n_in_j(q_t+v_kq_{,k})
\notag\\&\quad
+[\xi\lambda_1 (1-q^2)/\Gamma-\zeta]\Delta\chi q (n_i n_j-\delta_{ij}/2) 
.
\label{new:extrastressl}
\end{align}
The Leslie constants $\alpha_i$ and the parameters of Beris-Edwards model are
related by (see (2.10)-(2.15) in~\cite{DOY01})
\begin{subequations}\label{alqconv}
\begin{align}
\alpha_1(q)&=-\frac{2}{3}q^2(3+4q-4q^2)\xi^2/\Gamma,\\
\alpha_2(q)&=\left\{-\frac{1}{3}q(2+q)\xi-q^2\right\}/\Gamma,\\
\alpha_3(q)&=\left\{-\frac{1}{3}q(2+q)\xi+q^2\right\}/\Gamma,\\
\alpha_4(q)&=\frac{4}{9}(1-q)^2\xi^2/\Gamma+{2\mu},\\
\alpha_5(q)&=\left\{\frac{1}{3}q(4-q)\xi^2+\frac{1}{3}q(2+q)\xi\right\}/\Gamma,\\
\alpha_6(q)&=\left\{\frac{1}{3}q(4-q)\xi^2-\frac{1}{3}q(2+q)\xi\right\}/\Gamma.
\end{align}
\end{subequations}

We conclude that under representation \rf{Qrep} the  model
\rf{cont1}--\rf{mol1} turns into four equations
\begin{subequations}
\label{ErMod}
\begin{align}
0&=\pai v_i, \label{cont2}\\
0&=-\pai p -L_1\pj \left(\frac{3}{4}q_{,i}q_{,j}+2q^2n_{k,i}n_{k,j}\right)
+ \pj\tilde T_{ij}
,
\label{mom2} \\
L_1\Gamma[2qn_{i,kk}+4q_{,k}n_{i,k}]&=2qN_i{-2qL_1\Gamma|n_{j,k}|^2n_i}
\notag\\&\quad
-\frac{2}{3}(q+2)\xi[e_{ji}n_j-e_{lk}n_ln_kn_i] ,
\label{mol2n}\\
q_t+v_kq_{,k}&=\frac{2}{3}(q+2)\xi e_{ji}n_jn_i + L_1{\Gamma}q_{,kk}{-4qL_1\Gamma|n_{j,k}|^2}
\notag\\&\quad
+\Gamma \left(a^2q-\frac{c^2q^3}{2}\right)+\lambda_1\Delta\chi q,
\label{mol2q}		
\end{align}  
\end{subequations}
where $\tilde T_{ij}$ is given by \rf{new:extrastressl}. 

\paragraph{Boundary conditions at the substrate.}
Using \eqref{Qrep} in \eqref{Q:bcsubs}, the boundary conditions at $x_3=0$ become
\begin{subequations}
\label{ErMod_BC1}
\begin{align}
v_1&=0, \quad v_3=0,\\
n_3 &=\cos\theta_1,\\
q &= q_1.
\end{align}
\end{subequations}

\paragraph{Boundary conditions at the free interface.}

The condition \eqref{Q:Ts1}
now takes the form
\be\label{eqn:tsij_mod}
T_{ij}^s=g_0(\delta_{ij}-\nu_i\nu_j).
\ee
Projecting this condition onto the normal and tangential directions at the interface gives 
\begin{subequations}
\label{ErMod_BC2}
\begin{align}
\nu_i T_{ij}\nu_j &= -g_0 \pai\nu_i\label{e:bcn-flow}\\
\nu_i T_{ij} t_j &=0. \label{e:bct-flow}
\end{align}
The remaining conditions at the free interface $x_3=\eta(x_1,t)$ are
\begin{align}
\eta_t &= v_3-v_1 \partial_1\eta,\\
n&=R(\theta_2)\nu,\\
q&=q_2.
\end{align}
\end{subequations}

\section{Derivation of thin-film models}\label{sec:lub-ericksen}

\subsection{Thin-film model for active the Erickson theory}
We now non-dimensionalize this model using length scales $L$ for $x_1$ and
$\eps L$ for $x_3$, where $L$ denotes the characteristic lateral extend of the
cell and $\eps L$ denote its height. Hence, $\eps$ is the ratio between the two
length sales and in a thin-film settimg assumed to be small. We denote 
\begin{equation}
\lb{sc_b}
\begin{array}{cccccccccccccccc}
x_3&=&\eps L\bar x_3,&\quad
x_1&=&L\bar x_1,&\quad
\eta&=&\eps L\bar\eta,&\quad
&&&&\\
v_1&=&U\bar v_1,&\quad
v_3&=&V\bar v_3,&\quad
t&=&({L}/{U})\,\bar t,&\quad
&&&&\\
p&=&p_0+P \bar p,&\quad
h&=&\wscale \bar h,&\quad
&&&&
\end{array}
\end{equation}
where ${\cal E}$ and  $P$ are 
are defined as
\be
\wscale=\frac{L_1 }{\eps^2 L^2}\label{new:Gamma},
\ee
\be
P= \frac{\mu U}{\eps^2 L}\label{new:P}.
\ee
The order parameter $q$ and the the director field $n$ are dimensionless and do not need to be scaled. In the normal stress condition at the free surface, balancing the pressure with surface tension requires 
\be
P=\frac{\eps\,g_0}{L}.\label{new:P2}
\ee
Together with \eqref{new:P} this means 
\be\label{new:eps}
\eps^3=\frac{\mu U}{g_0}.
\ee
Further scalings are obtained as
\begin{subequations}
\begin{align}
N_i&=\frac{U}{\varepsilon L}\bar N_i,\\
e_{11}&= \frac{U}{L} \bar e_{11},
\quad
e_{13}=\frac{U}{\varepsilon L} \bar e_{13},
\quad
e_{31}=\frac{U}{\varepsilon L} \bar e_{31},
\quad
e_{33}=\frac{U}{L} \bar e_{33},\\
\omega_{13}&=\frac{U}{\varepsilon L} \bar \omega_{13},
\quad
\omega_{31}=\frac{U}{\varepsilon L} \bar \omega_{31},
\\
\alpha_i&=\mu\bar\alpha_i
\quad
\Gamma=\bar\Gamma/\mu,\quad
a^2={\cal E}\bar a^2,
\quad
c^2={\cal E}\bar c^2,\\
\tilde T_{ij}&= \frac{\mu U}{\varepsilon L}\bar{\tilde T}_{ij},\quad
\left[T_{11}^E,\, T_{13}^E,\, T_{31}^E,\, T_{33}^E\right]= \frac{\mu U}{\varepsilon L}
\left[\eps^2 T_{11}^E,\, \eps T_{13}^E,\, \eps T_{31}^E,\, T_{33}^E\right],
\\
\bar L_1 &= \frac{L_1}{\varepsilon \mu U L}, 
\quad \bar\zeta\Delta\bar\chi=\frac{\Gamma L}{U}\zeta\Delta\chi,
\quad \bar\lambda_1\Delta\bar\chi=\frac{L}{U}\lambda_1\Delta\chi.
\lb{L_chi_scaling}
\end{align}
\lb{sc_f}
\end{subequations}

Retaining only the leading order terms in $\lubp$ in the rescaled system \rf{ErMod}, given in appendix A, and assuming the weak
elasticity limit, by which we mean that as we introduce the thin-film approximation $\eps\to0$, we assume $\bar{L}_1 = O(1)$
and $\Delta\bar{\chi} = O(\eps^{-1})$, the leading order system in the bulk becomes 
\begin{subequations}
\begin{align}
0&=v_{1,1}+v_{3,3}, \label{be:cont_lo}\\
0&=-p_{,1}+\frac{1}{2}(v_{1,3}f_A(n_1,n_3))_{,3}
\notag\\&\qquad
+{\frac{\eps}{\Gamma}\Delta\chi\left[\left(\xi\lambda_1 (1-q^2)-\zeta\right) q n_1 n_3 \right]_{,3}},
		\label{be:mom1_lo} \\
0&=-p_{,3},\label{be:mom2_lo} \\
L_1\Gamma[2qn_{1,33}+4q_{,3}n_{1,3}]&
=-2q\frac{v_{1,3}n_3}{2}-2qL_1\Gamma \left(|n_{1,3}|^2+|n_{3,3}|^2\right)n_1
\nonumber\\&\qquad
-\frac{2v_{1,3}}{3}(q+2)\xi[n_3/2-n_1^2n_3], 
\label{be:mol_n1_lo}\\
L_1\Gamma[2qn_{3,33}+4q_{,3}n_{3,3}]&
=2q\frac{v_{1,3}n_1}{2}-2qL_1\Gamma \left(|n_{1,3}|^2+|n_{3,3}|^2\right)n_3
\nonumber\\&\qquad
-\frac{2v_{1,3}}{3}(q+2)\xi[n_1/2-n_1n_3^2], 
\label{be:mol_n2_lo}\\
-\frac{2}{3}(q+2)\xi v_{1,3}n_1n_3&=
-4qL_1\Gamma  \left(|n_{1,3}|^2 +  |n_{3,3}|^2  \right) 
+
L_1\Gamma q_{,33}
\nonumber\\&\qquad
+\Gamma \left(a^2q-\frac{c^2q^3}{2}\right)
{+\eps\lambda_1\Delta\chi q},
\label{be:mol_q_lo}		
\end{align}  
\end{subequations}
where here and below for convenience we have skipped the overbars everywhere and in the horizontal momentum equation we introduced the notation
\begin{align}
f_A(n_1,n_3)&\equiv 2\alpha_1(n_1n_3)^2+(\alpha_5-\alpha_2)n_3^2+(\alpha_3+\alpha_6)n_1^2+\alpha_4.\end{align}
The leading order system for the boundary conditions at $z=0$ is given by 
\begin{subequations}
\begin{align}
v_1&=0, \quad v_3=0,\\
n_3 &=\cos\theta_1,\\
q &= q_1,
\end{align}
\lb{BC_z0}
\end{subequations}
and at the free surface, $x_3=\eta(x_1,t)$ by (cf. \rf{NSCond} in Appendix A).
\begin{subequations}
\begin{align}
\eta_t &= v_3-v_1 \partial_1\eta,\\
-p&= \eta_{,11} ,\\
\frac{1}{2}v_{1,3}f_B(n_1,n_3) &= -{
\frac{\eps}{\Gamma}\Delta\chi\left(\xi\lambda_1 (1-q^2)-\zeta\right) q n_3 n_1 },\\
n_3&=\cos\theta_2,\\
q&=q_2,
\lb{BC_znu}
\end{align}
\end{subequations}
with a given function $q_2(x,t)$, and where we define
\begin{align}
f_B(n_1,n_3)&\equiv2\alpha_1(n_1n_3)^2+(\alpha_6-\alpha_3)n_3^2+(\alpha_2+\alpha_5)n_1^2+\alpha_4.
\end{align}
Next, similar to~\cite{Lin2013} we rewrite equations \rf{be:cont_lo}--\rf{be:mol_q_lo} in the bulk in terms of the director angle $\theta$ using the representation 
\be
n_1=\sin\theta,\quad n_3=\cos\theta.
\lb{n_Rep}
\ee
First, let us multiply equations \rf{be:mol_n1_lo} and \rf{be:mol_n2_lo} by $-n_3$ and $n_1$, respectively, and then sum up them. The resulting equation has the form:
\begin{equation}
L_1\Gamma[-2qn_{1,33}n_3-4q_{,3}n_{1,3}n_3+2qn_{3,33}n_1+4q_3n_{3,3}n_1]=qv_{1,3}+\frac{v_{1,3}}{3}(q+2)\xi\left[n_3^2-n_1^2\right].
\end{equation}
Using \rf{n_Rep} and definitions 
\be
\gamma_1=\alpha_3 - \alpha_2,\quad \gamma_2=\alpha_2+\alpha_3=\alpha_6-\alpha_5
\label{new:Onsager}
\ee
and \rf{alqconv} the 
latter equation can be reduced to \eqref{be:mol_n_lo_theta} which we include
with the other equations of the system \eqref{be:cont_lo}--\rf{be:mol_q_lo},
rewritten in terms of $\theta$, giving
\begin{align}
0&=v_{1,1}+v_{3,3}, \label{be:cont_lo_theta}\\
0&=-p_{,1}+\frac{1}{2}(v_{1,3}f_A(\theta))_{,3}
\notag\\ & \quad 
+{
\frac{\eps}{2\Gamma}\Delta\chi\left[\left(\xi\lambda_1 (1-q^2)-\zeta\right) q \sin(2\theta) \right]_{,3}},
		\label{be:mom1_lo_theta} \\
0&=-p_{,3},\label{be:mom2_lo_theta} \\
L_1q\left[2q\theta_{,33}+4q_{,3}\theta_{,3}\right]&
=-\frac{v_{1,3}}{2}[\gamma_1-\gamma_2\cos(2\theta)]
\label{be:mol_n_lo_theta}\\
-\frac{1}{3}(q+2)\xi v_{1,3}\sin(2\theta)&=-4qL_1\Gamma|\theta_{,3}|^2
\notag\\ & \quad 
+L_1\Gamma q_{,33}+\Gamma(a^2q-\frac{c^2q^3}{2})
{ +\eps\lambda_1\Delta\chi q},
\label{be:mol_q_lo_theta}		
\end{align}  
where we define
\begin{align}
f_A(\theta)&=(\alpha_1/2)\sin^2(2\theta)+(\alpha_5-\alpha_2)\cos^2\theta+(\alpha_3+\alpha_6)\sin^2\theta+\alpha_4.
\end{align}
The leading order system for the boundary conditions at $z=0$ is given by 
\begin{subequations}
\begin{align}
v_1&=0, \quad v_3=0,   \label{eqntheta:vbot}\\
\theta&=\theta_1,\\
q &= q_1,
\end{align}
\lb{eqntheta:BC_z0}
\end{subequations}
and at the free surface, $x_3=\eta(x_1,t)$ by
\begin{subequations}
\begin{align}
\eta_t &= v_3-v_1 \partial_1\eta, \label{eqntheta:kincon}\\
-p&= \eta_{,11} ,\label{eqntheta:eta11}\\
\frac{1}{2}v_{1,3}f_B(\theta) &=
-{
\frac{\eps}{2\Gamma}\Delta\chi\left(\xi\lambda_1 (1-q^2)-\zeta\right) q
\sin(2\theta) }
,\label{eqntheta:ts}\\
\theta&=\theta_2,\\
q&=q_2,
\lb{eqntheta:BC_znu}
\end{align}
\end{subequations}
where we define
\begin{align}
f_B(\theta)&=(\alpha_1/2)\sin^2(2\theta)+(\alpha_6-\alpha_3)\cos^2\theta+(\alpha_2+\alpha_5)\sin^2\theta+\alpha_4.
\end{align}

We now integrate these equations. First, the combination of 
\eqref{be:cont_lo_theta} and \eqref{eqntheta:kincon}
gives
\bes
\eta_t(x_1,t)=-\partial_1 \int_0^\eta v_1(x_1,x_3,t)dx_3.
\ees
which is in fact exact i.e.\ also valid for the full governing equations.
From \eqref{be:mom2_lo_theta}, \eqref{eqntheta:eta11},
\eqref{be:mom1_lo_theta}, \eqref{eqntheta:ts}, 
we get
\begin{align}
p&=-\eta_{,11},\\
 f_A(q,\,\theta) v_{1,3}
&=2\eta_{,111}(\eta-x_3)
-{
 \frac{\eps}{\Gamma}\Delta\chi\left(\xi\lambda_1 (1-q^2)-\zeta\right) 
q \sin(2\theta)}
\notag\\&\quad
+{ \frac{\eps}{\Gamma}\Delta\chi\left(\xi\lambda_1 (1-q_2^2)-\zeta\right) q_2 
 \sin(2\theta_2)}
\notag\\&\quad
-{
\frac{\eps}{\Gamma}\Delta\chi
\frac{f_A(q_2,\theta_2)}{f_B(q_2,\theta_2)}
\left(\xi\lambda_1 (1-q_2^2)-\zeta\right) q_2 \sin(2\theta_2) }
\notag\\
&=2\eta_{,111}(\eta-x_3)
-{
\frac{\eps}{\Gamma}\Delta\chi\left(\xi\lambda_1 (1-q^2)-\zeta\right) 
q \sin(2\theta)}
\notag\\&\quad
+{
\frac{\eps}{\Gamma}\Delta\chi
\frac{\gamma_1-\gamma_2 \cos(2\theta_2)}{f_B(q_2,\theta_2)}
\left(\xi\lambda_1 (1-q_2^2)-\zeta\right) q_2 \sin(2\theta_2) }
,\label{eqn:li1:f2v13}
\end{align}
provided $f_B(q_2,\theta_2)\neq0$.

{Therefore, the last three equations combined together result in a closed\\ {\it lubrication system}:}
\begin{subequations}
\lb{Lub_mod}
\begin{align}
\label{eqn:etat}
\eta_t(x_1,t)&=-\partial_1 \int_0^\eta v_1(x_1,x_3,t)dx_3,\\
v_{1,3}&=\frac{2\eta_{,111}}{f_A(q,\theta)}(\eta-x_3)
\notag\\&\quad
-\frac{\eps\Delta\chi}{{\Gamma}f_A(q,\theta)}
\Big[
{
\left(\xi\lambda_1 (1-q^2)-\zeta\right) 
q \sin(2\theta)}
\notag\\&\qquad\qquad\qquad
\left.
-{
\frac{\gamma_1-\gamma_2 \cos(2\theta_2)}{f_B(q_2,\theta_2)}
\left(\xi\lambda_1 (1-q_2^2)-\zeta\right) q_2 \sin(2\theta_2) }
\right],
\label{eqn:li1:v13}\\
\left(q^2\theta_{,3}\right)_{,3}&
=-\frac{1}{4L_1}\left(\gamma_1-\gamma_2\cos(2\theta)\right)v_{1,3},
\label{eqn:li1:mol_n_lo_theta}\\
q_{,33}
&=
4q(\theta_{,3})^2
-\frac{\xi(q+2)}{3L_1\Gamma} \sin(2\theta)v_{1,3}
-\frac q{L_1}\left(a^2-\frac{c^2q^2}{2}\right)
-{\frac{\eps\lambda_1\Delta\chi }{L_1\Gamma} q}.
\label{eqn:li1:mol_q_lo_theta}	
\end{align}

Notice that if $f_A(q(x_1,x_3,t),\,\theta(q(x_1,x_3,t))\neq 0$ 
for all $x_1$, $x_3$ and $t$,
we can solve \eqref{eqn:li1:f2v13} for $v_{1,3}$ and use the result in 
\eqref{be:mol_n_lo_theta} and \eqref{be:mol_q_lo_theta} to eliminate 
$v_{1,3}$, thus decoupling the system for $\theta$ and $q$ from the velocity field.
Because of the size of the resulting equations, we have not done this here.

\begin{align}
v_1&=0 \qquad \text{at } x_3=0,\\
\lb{theta_BC}
\theta&=\theta_1 \qquad \text{at } x_3=0,\\
q&=q_1 \qquad \text{at } x_3=0,\\
\theta&=\theta_2 \qquad \text{at } x_3=\eta,\\
\lb{q_BC}
q&=q_2 \qquad \text{at } x_3=\eta.
\end{align}
\end{subequations}

\subsection{Thin-film model for the active Leslie-Erickson-Parodi theory}
If we use the Leslie-Erickson-Parodi theory with correspnding active terms as a model for the active liquid crystal \cite{KJJPS05,Juelicher2007,VJP05,Sankararaman2009} and nondimensionalise as before we derive in appendix B the following coupled system for the leading order thin-film approximation 
\begin{subequations}
\begin{align}
\partial_t \eta&= -\pone\,\int^{\eta}_0  v_1\,dx_3,  \label{kin}\\
0&=\eta_{,111} (x_3-\eta)+\frac{1}{2}v_{1,3}f_A(\theta)
+\frac{\zeta^{ELP}\Delta\chi^{ELP}}{2}
     \frac{f_A(\theta_2)}{f_B(\theta_2)}\sin(2\theta_2)
\notag \\&\quad
+\frac{\zeta^{ELP}\Delta\chi^{ELP}}{2}\left(\sin(2\theta)-\sin(2\theta_2)\right)
,
\\
2K \theta_{,33} &= -(\gamma_1- \gamma_2 \cos(2\theta)) v_{1,3},
\end{align}
\lb{LEP_LM}
\end{subequations}
with the boundary conditions at $x_3=0$ given by
\begin{subequations}
\begin{align}
\quad v_3&=0,   \\
\theta&=\theta_1,
\end{align}
\end{subequations}
and at the free surface, $x_3=\eta(x_1,t)$, 
\begin{align}
\theta&=\theta_2.
\end{align}

Formal comparison of \rf{LEP_LM} with equations \rf{eqn:etat}-\rf{eqn:li1:mol_n_lo_theta}
considered with $q=q_1=q_2=\const$ provides the following relations between
the active and elastic parameters in the Eriksen-Leslie-Parodi and Ericksen thin-film models:
\be
\lambda^{ELP}_1=\lambda_1,\ \zeta^{ELP}=(\zeta-\xi\lambda_1(1-q^2))q,\ K=2L_1q^2,\ \Delta\chi^{ELP}=\frac{\eps}{\Gamma}\Delta\chi.
\ee

At the same time, in absence of the active terms $\Delta\chi^{ELP}=0$ our model \rf{LEP_LM} can be shown to coincide with the (passive) thin-film model derived in~\cite{Lin2013} for the weak elasticity regime, cf. system (A17)-(A20) there. Note that the special anchoring boundary conditions $\theta_1=\pi/2$ and $\theta_2=0$  were considered  in Lin et al. \cite{Lin2013}.

\section{Impact of activity terms}
At this point, further reductions of the thin-film model \rf{Lub_mod} or \rf{LEP_LM} are not, in general, possible without additional assumptions, since the remaining equations cannot be easily integrated with respect to $x_3$.
We will instead look at two special cases of the more general $Q$-tensor system \rf{Lub_mod}: one, where the interface is flat ($\eta=1$)
and the other where the misalignment of the director at the substrate and the interface
is small, $|\theta_2-\theta_1|\ll 1$. 
\subsection{Flat film}
\paragraph{Passive case}
We first consider the case, where $\eta=\const$ is any positive constant. This yields $v_1=0$ and 
\begin{subequations}
\begin{align}
\lb{theta_r}
q^2 \theta_{,3} &= c_1,\\
\lb{theta_q}
q_{,33} &= 4q  (\theta_{,3})^2 - \frac1{L_1}\left(a^2q-\frac{c^2q^3}{2}\right) 
\end{align}
\end{subequations}
Under the additional assumption that $q_1=q_2\equiv q_0$ and that $q$ remains constant we obtain the solution
\begin{equation}
\lb{theta_rt}
\theta=(\theta_2-\theta_1) x_3 +\theta_1,
\qquad
q=q_0=\left[\frac{2a^2}{c^2}-\frac{8 L_1}{c^2}\left(\theta_2-\theta_1\right)^2\right]^{1/2}.
\end{equation}
We note that a similar solution for the director angle $\theta$ and for $q$ has been found for the case of channel flow in \cite{Batista2015}.

Alternatively, one can also substitute \rf{theta_r} into \rf{theta_q} to
obtain one ODE for $q$:
\be
q_{,33}=\frac{4c_1^2}{q}-\frac{1}{L_1}(a^2q-\frac{c^2q^3}{2}).
\lb{eq_q1}
\ee
Multiplying the last equation by $q_{,3}$ and integrating in $x_3$ one obtains
\bes
\frac{1}{2}q_{,3}^2=4c_1^2\log(q)-\frac{1}{L_1}(\frac{a^2}{2}q^2-\frac{c^2q^4}{8})+c_2,
\ees
where we have assumed that $q\not=\const$ and
\be
c_2=\left(-4c_1^2\log(q_1)+\frac{1}{L_1}(\frac{a^2}{2}q_1^2-\frac{c^2q_1^4}{8})+\frac{1}{2}q_{,3}^2\right)\Big|_{x_3=0}.
\lb{c_2}
\ee
The last ODE is separable and can be integrated as
\be
x_3=\int_{q_1}^q\frac{ds}{\sqrt{8c_1^2\log(s)-\frac{1}{L_1}(a^2s^2-\frac{c^2s^4}{4})+2c_2}},
\lb{q_sol}
\ee
where we have assumed that $q_2>q_1$.
Correspondingly, using \rf{theta_r} one finds
\be
\theta(x_3)-\theta_1=\int_0^{x_3}\frac{c_1}{q^2(x_3)}\,dx_3
=c_1\int_{q_1}^q\frac{ds}{s^2\sqrt{8c_1^2\log(s)-\frac{1}{L_1}(a^2s^2-\frac{c^2s^4}{4})+2c_2}}.
\lb{theta_sol}
\ee
In the last expression the constants $c_1$ and $c_2$ are determined by the boundary condition
for $\theta$ at $x_3=\eta=\const$:
\bea
\theta_2-\theta_1&=&c_1\int_{q_1}^{q_2}\frac{ds}{s^2\sqrt{8c_1^2\log(s)-\frac{1}{L_1}(a^2s^2-\frac{c^2s^4}{4})+2c_2}},\nonumber\\
\eta&=&\int_{q_1}^{q_2}\frac{ds}{\sqrt{8c_1^2\log(s)-\frac{1}{L_1}(a^2s^2-\frac{c^2s^4}{4})+2c_2}}.
\lb{CompCond}
\eea
The compatibility conditions \rf{CompCond} do not have always solutions. For example,
if $\theta_2=\theta_1$ and $q_1$ is large enough the denominator in
\rf{CompCond} is non-negative for all $q\ge q_1$ and the first integral in \rf{CompCond} can not be zero. Therefore, in this case one has only the trivial solution \rf{theta_rt}.

\paragraph{Active flat film}
The compatibility condition of $\eta(x_1,x_3,t)=\eta=\const$ with  \rf{eqn:etat} implies that $v_1,\,q,\,\theta$ are functions of $x_3$ only.
By that, system \eqref{eqn:li1:v13}-\rf{eqn:li1:mol_q_lo_theta} reduces to 
\begin{subequations}
\label{flat_active}
\begin{align}
v_{1,3}&=-\frac{\eps\Delta\chi}{\Gamma f_A(q,\theta)}
\Big[
\left(\xi\lambda_1 (1-q^2)-\zeta\right) q \sin(2\theta)\Big.\\
&\hspace*{3cm}  \Big. -\frac{\gamma_1-\gamma_2 \cos(2\theta_2)}{f_B(q_2,\theta_2)}
\left(\xi\lambda_1 (1-q_2^2)-\zeta\right) q_2 \sin(2\theta_2)\Big]
\label{eqn:li1:v14}\\
\left(q^2\theta_{,3}\right)_{,3}&
=-\frac{1}{4L_1}\left(\gamma_1-\gamma_2\cos(2\theta)\right)v_{1,3},
\label{eqn:li1:mol_n_lo_theta_1}\\
q_{,33}
&=
4q(\theta_{,3})^2
-\frac{\xi(q+2)}{3L_1\Gamma} \sin(2\theta)v_{1,3}
-\frac q{L_1}\left(a^2-\frac{c^2q^2}{2}\right)
-\frac{\eps\lambda_1\Delta\chi }{L_1\Gamma} q,
\label{eqn:li1:mol_q_lo_theta_1}	
\end{align}
\end{subequations}
which further reduce to two coupled ODEs for $\theta(x_3)$ and $q(x_3)$ by eliminating $v_{1,3}$. The latter ODEs can be effectively integrated numerically.

Note, that in absence of the active terms ($\lambda_1=0$ or $\zeta=0$) the nontrivial solution to the system \rf{flat_active} is given by \rf{q_sol}-\rf{theta_sol} combined with $v_1=0$ and it exists only when the compatibility conditions \rf{CompCond}  on the boundary data \rf{theta_BC}-\rf{q_BC} are satisfied. Given  $q_2>q_1$, such that the square root in the denominator of  \rf{q_sol} is real for all $q\in (q_1,\,q_2)$,
by taking $\eta$ and $\theta_2-\theta_1$ sufficiently large, one can realize this  passive solution. Moreover, 
also for small active terms with $\Delta\chi\ll 1$
this non-homogeneous solution to the system \rf{flat_active} continously persists and by \rf{eqn:li1:v14} exibits the non-homogeneous flow $v_1(x_3)$ with $|v_1|\ll 1$. This effect of inducing a non-zero flow in a channel geometry, when the thickness of the latter $\eta$ becomes sufficiently large, has been observed in~\cite{VJP05,VJP06} for the polar Leslie-Ericksen-Parodi based models.

Finally, note that when active terms are present in \rf{flat_active} there is no analogous solution to \rf{theta_rt}. One can show that the ansatz 
\rf{theta_rt} does not satisfy equations \rf{eqn:li1:v14}-\rf{eqn:li1:mol_n_lo_theta_1}, unless $q_0=0$.

\subsection{Film with small angle change in the director boundary condition} 
Another special case, where it is possible to discuss analytical solutions is obtained if the difference in the director angle is small.
\paragraph{Passive case}
Assuming $|\theta_2-\theta_1|\ll 1$, then to leading
order $\theta=\theta_2=\theta_1$ is constant and  
\rf{eqn:li1:mol_n_lo_theta} imples $v_{1,3}=v_1=0$ and $\eta=\const$. As a result the whole dynamics reduces to \rf{eqn:li1:mol_q_lo_theta}, which can be further reduced to \rf{eq_q1} with $c_1=0$. Then the corresponding
solution is given by
\be
x_3=\frac{1}{\sqrt{2}}\int_{q_1}^q\frac{ds}{-\frac{1}{L_1}(\frac{a^2}{2}s^2-\frac{c^2s^4}{8})+c_2}.
\lb{const_theta}
\ee
{The compatibility conditions \rf{CompCond} reduce to 
\be
\eta=\int_{q_1}^q\frac{ds}{\sqrt{8c_1^2\log(s)-\frac{1}{L_1}(a^2s^2-\frac{c^2s^4}{4})+2c_2}}
\lb{CompCond_r}
\ee }

We note that the solution \rf{const_theta} with $\theta=\theta_1=\const$ to the system \rf{eqn:etat},  \eqref{eqn:li1:v13}-\rf{eqn:li1:mol_q_lo_theta}
does not exist in the case when active terms are present ($\lambda_1\not=0$ or $\zeta\not=0$), since in that case \eqref{eqn:li1:v13}-\rf{eqn:li1:mol_n_lo_theta} are not satisfied.

\paragraph{Active film}
Another way to initiate a nontrivial dynamics in the case  $\theta=\theta_2=\theta_1=\const$
is to assume
\be
\gamma_1-\gamma_2\cos(2\theta)=0.
\lb{gamma_cos}
\ee
This would imply that \rf{eqn:li1:mol_n_lo_theta}  is satisfied and $q=q_1=q_2=\const$. Furthermore, \eqref{eqn:li1:v13} can be integrated and introduced
into \eqref{eqn:etat} yields a new modified thin-film equation
\begin{equation}
\eta_t=- \frac{2}{3f_A(q_1,\,\theta_1)}\partial_1 \left[\eta^3\eta_{,111}\right] +C(q_1,\,\theta_1)\partial_1(\eta^2).
\lb{lub_eq}
\end{equation}
Note that in this case, besides the trivial isotropic solution $q_1=0$, only special values of $q_1$ and $\theta_1$ are allowed. These have to be compatible with both, 
equation \rf{gamma_cos} and the algebraic relation 
\be
0=\left(a^2-\frac{c^2q^2}{2}\right)+\frac{\eps\lambda_1\Delta\chi}{\Gamma}.
\lb{q_lambda_1}
\ee
which arises from \rf{eqn:li1:mol_q_lo_theta}.

For given activity $\lambda_1\Delta\chi\in\mR$, solutions for $q_1$ and $\theta_1$ can be obtained from \rf{gamma_cos} and \rf{q_lambda_1} as
\be
\cos(2\theta)=-\frac{3}{\xi}+\frac{6}{(2+q)\xi},\quad q^2=\frac{2a^2}{c^2}+\frac{2\eps\lambda_1\Delta\chi}{c^2\Gamma}.
\lb{const_theta_AT}
\ee
Finally, note that solution \rf{const_theta_AT} to system \rf{Lub_mod} does not always exists. In particular, it does not exists for $\xi=0$,
i.e. when liquid crystal molecules align perfectly with the hydrodynamic flow. In absence of active terms ($\lambda_1=0$ or $\zeta=0$)
one has $C(q_1,\,\theta_1)=0$ in \rf{lub_eq}, and therefore the hydrodynamic flow decouples from the nematics via the rescaling of time by $f_A(q_1,\,\theta_1)$.

\section{Discussion and outlook}\label{sec:conc}   
In this article we presented a systematic asymptotic derivation of the thin-film model  given by the system \rf{Lub_mod} from the free-boundary problem for the Beris-Edwards model to describe the evolution and flow structure of an active nematic liquid crystal. We also showed that the new thin-film model formally reduces to the polar one based on Leslie-Ericksen-Parodi theory \rf{LEP_LM}, when the scalar order parameter $q$ is homogeneous, which  in the passive case coincides with the model  derived previously in~\cite{Lin2013}.  In the active case
our analytical  solution to  \rf{Lub_mod} demonstrates nonzero flow that can be spontaneously initiated from the homogeneous one by increasing the film thickness, as previously observed in~\cite{VJP05,VJP06}. 

We now point to some further applications as well as extensions of our results.
The derivation of the coupled model \rf{Lub_mod} 
starting from the Ericksen type model \rf{ErMod} considered with boundary
conditions \rf{ErMod_BC1}--\rf{ErMod_BC2} has been conveyed under
the assumption of continuity of the director field $n$ in the film bulk and at the free surface. We note that these models are capable to describe solutions having point defects of integer degree $k$ with $k\in\mZ$. Two typical examples of defects with degree $-1$ in the film bulk and of degree $1$ at the film contact line are presented in Fig. 2. 
\begin{figure}
\centering
\includegraphics*[width=0.58\textwidth]{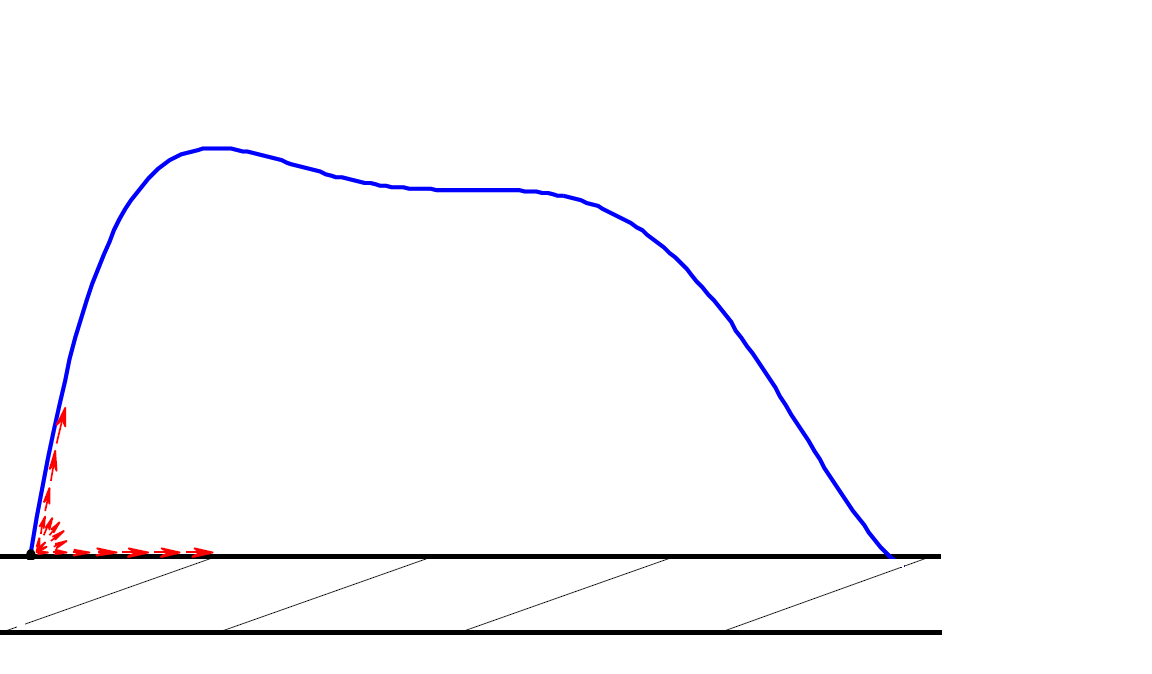}
\includegraphics*[width=0.38\textwidth]{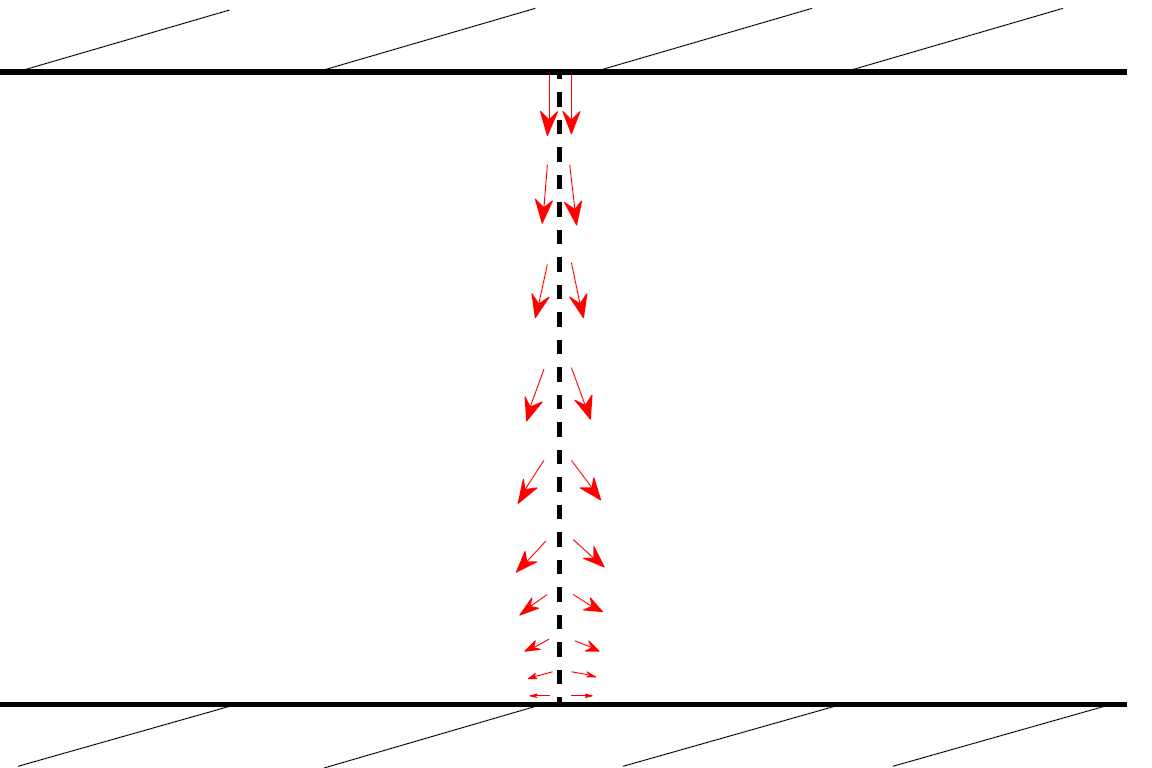}
\caption{Examples of defects of degree $1$ at the contact line (left) and a wall defect (right), which mathematically can be described
by the point defect of degree $-1$ located at the intersection of the wall (depicted by dashed line) and the substrate. The direction of the
director filed in the neighborhood of the defect is shown by red arrows. The magnitude of scalar order parameter $q$ is represented by the arrow size.
The derived $Q$-tensor lubrication model  \rf{Lub_mod} can smoothly resolve such integer point defects through continous reduction of the scalar parameter field to $q=0$ when approaching
the defect.} 
\label{Fig1}
\end{figure}
One observes that when approaching the defect points the magnitude of the scalar order parameter $q$ goes to zero and by that preserving the continuity of the full $Q$-tensor field \rf{Qrep}. 
We should also point out that, being derived in the weak elasticity regime (cf. scaling for $L_1$ in \rf{L_chi_scaling}) and under the large pressure scaling \rf{new:P}, the model \rf{eqn:etat},\rf{Lub_mod}  allows for $O(1)$ variation of the director field $n$ along the vertical $x_3$ direction of the film. This is the case, for example, in the wall defect of degree $-1$ in a confined flat film presented in Fig. 1 (right plot), where the director angle $\theta$ changes from $0$ to $\pi/2$ along the vertical film direction. Such wall defects were observed before in experiments on thin passive nematic films~\cite{LP95}.

However, the Ericksen model may not always resolve defects of rational degree $k+1/2,\,k\in\mZ$, because the latter exibit special disclination lines along which $n$ changes to $-n$~\cite{FRSZ16,KRSZ16, KA15}. Nevertheless, such defects can be described by the lubrication model \rf{Lub_mod} if the special condition
\be
\displaystyle\lim_{(x_1,\,x_3)\rightarrow(x_1^*,\,x_3^*)}|\sqrt{q(x_1,\,x_3)}\nabla n(x_1,\,x_3)|<\infty
\label{Def_Reg_Cond}
\ee
is fulfilled, where $(x_1^*,\,x_3^*)$ is an instant defect location. It is easy to check that \rf{Def_Reg_Cond} ensures then that the associated local Ericksen elastic energy and the corresponding terms involving  director gradients in \rf{eqn:li1:mol_n_lo_theta}--\rf{eqn:li1:mol_q_lo_theta} are kept finite.

In the Ericksen theory the defects are defined as singular points where the
scalar order parameter $q=0$~\cite{Za12}. We note, that besides the 
singularities of the director field, considered in the previous point,
the solutions to \rf{Lub_mod} may exibit singular lines along
which $q=0$ but $n$ ist still continous. These lines have special physical
meaning, because $Q$-tensor in \rf{Qrep} is zero and, therefore, the nematic field is isotropic along them.
A typical example of such a line is given by the middle line of the channel at the right imbedded $Q$-tensor plot in Fig. 2 of~\cite{Batista2015}.

In the future, we plan to investigate system \rf{Lub_mod} numerically, where the nematic part of the model, equations \rf{eqn:li1:mol_n_lo_theta}-\rf{eqn:li1:mol_q_lo_theta}, is
given by two coupled ODEs in $x_3$ direction. 
A quasi one-dimensional
numerical scheme could then be developed by solving the hydrodynamic equations  \rf{eqn:etat}, \rf{eqn:li1:v13}
and the nematic ones \rf{eqn:li1:mol_n_lo_theta}-\rf{eqn:li1:mol_q_lo_theta} separately and iteratively. We expect complicated
solution patterns for \rf{Lub_mod} to arise
with dynamical formation, mutual interaction and anighilation of point
defects in the film bulk similar to ones observed recently in~\cite{Thampi2014b,ZZRP16}. As in the latter
works, it would be important to investigate the interrelation of these patterns
with the evolution of the liquid vorticity field~\cite{Thampi2016}.

Finally, we note that by imposing the constant scalar order parameter $q=q_2$ in \rf{con_bc} we neglected possible Marangoni effects at the free surface and, in particular, in condition \rf{e:bct-flow} for the tangential stress. This was motivated by the fact that under the balance \rf{new:eps}, that keeps the surface tension term at leading order in \rf{e:bcn-flow}, the equation for the Marangoni force (see e.g. formula (8) in~\cite{Rey99}) necessarily impies that $q$ should be constant at the film free surface. Nevertheless, by relaxing condition \rf{new:eps} and neglecting surface tension one would be able to derive a model analogous to \rf{Lub_mod} for pure Marangoni driven active nematic thin films.

\section*{Acknowledgement}
The authors would like to thank Jonathan Robbins, Valeriy Slastikov and Arghir Zarnescu for their valuable comments on the results of this article. GK would like to acknowledge support from Leverhulme grant RPG-2014-226.  GK gratefully acknowledges the hospitality of the Weierstrass-Institute Berlin and OCIAM at the University of Oxford during his research visits.

\begin{appendix}
\section{Thin-film approximation}

After application of scalings \rf{sc_b}--\rf{sc_f} to Ericksen system \rf{ErMod} the non-dimensional equations for the bulk, after skipping overbars everywhere, take the form:
\begin{align}
0&=v_{1,1}+v_{3,3}, \label{be:cont2nondim}\\
0&=-p_{,1}-\eps^3 T^E_{11,1}-\eps T^E_{13,3}+\eps \tilde T_{11,1}+\tilde T_{13,3}
\label{be:mom21} \\
0&=-p_{,3} -\eps^3T^E_{31,1}-\eps T^E_{33,3}+\eps^2\tilde T_{31,1}+ \eps\tilde T_{33,3}
\label{be:mom22} 
\end{align}  
\begin{gather}		
L_1\Gamma\left(2\eps^2 qn_{1,11}+2q n_{1,33}+4\eps^2 q_{,1}n_{1,1}+4q_{,3}n_{1,3}\right)
\notag\\ 
=2qN_1-2qL_1
\Gamma\left(\eps^2 |n_{1,1}|^2+|n_{1,3}|^2+\eps^2 |n_{3,1}|^2+|n_{3,3}|^2\right)n_1
\notag\\ 
-\frac{2}{3}(q+2)\xi\left(\eps e_{11}n_1+e_{31}n_3
  -\left(\eps e_{11}n_1n_1+e_{13}n_1n_3+e_{31}n_3n_1+\eps e_{33}n_3n_3\right)n_1
\right);
\label{be:mol2n1}
\\[1ex]
L_1\Gamma\left(2\eps^2 qn_{3,11}+2q n_{3,33}+4\eps^2 q_{,1}n_{3,1}+4q_{,3}n_{3,3}\right)
\notag\\ 
=2qN_3-2qL_1
\Gamma\left(\eps^2 |n_{1,1}|^2+|n_{1,3}|^2+\eps^2 |n_{3,1}|^2+|n_{3,3}|^2\right)n_3
\notag\\ 
-\frac{2}{3}(q+2)\xi\left(e_{13}n_1+\eps e_{33}n_3
  -\left(\eps e_{11}n_1n_1+e_{13}n_1n_3+e_{31}n_3n_1+\eps e_{33}n_3n_3\right)n_3
\right);
\label{be:mol2n2}
\\[1ex]
\eps \left(q_t+v_kq_{,k}\right)
-\frac{2}{3}(q+2)\xi
\left( 
\eps e_{11}n_1n_1+e_{13}n_1n_3+e_{31}n_3n_1+\eps e_{33}n_3n_3
\right)
\notag\\[1ex]
=
-4q L_1\Gamma\left(\eps^2 |n_{1,1}|^2+|n_{1,3}|^2+\eps^2 |n_{3,1}|^2+|n_{3,3}|^2\right)+
\notag\\ \quad
L_1\Gamma\left(\eps^2q_{,11}+q_{,33}\right)+
L_1\Gamma\left(a^2q-\frac{c^2q^3}{2}\right)
+\eps\lambda_1\Delta\chi q
,
\label{be:mol2q}		
\end{gather}  
where
\begin{subequations}
\begin{align}
T_{11}^E=&L_1\left(\frac{3}{4}|q_{,1}|^2+2q^2|n_{k,1}|^2\right),
\lb{TE_11}\\
T_{13}^E=T_{31}^E=&L_1\left(\frac{3}{4}q_{,1}q_{,3}+2q^2n_{k,1}n_{k,3}\right),
\lb{TE_13}\\
T_{33}^E=&L_1\left(\frac{3}{4}|q_{,3}|^2+2q^2|n_{k,3}|^2\right),
\lb{TE_33}\\
\tilde{T}_{11}=&
\alpha_1 \left(
\eps n_1 n_1 e_{11} n_1 n_1 +
n_1 n_3 e_{13} n_1 n_1 +
n_3 n_1 e_{31} n_1 n_1 +
\eps n_3 n_3 e_{33} n_1 n_1 \right)
\notag\\  &\quad
+ \alpha_2 N_1n_1 + \alpha_3N_1n_1 
+\alpha_4 \eps e_{11}
\notag \\ &\quad 
+\alpha_5 \left(\eps e_{11}n_1n_1+e_{13}n_3n_1\right) 
+ \alpha_6 \left(\eps e_{11}n_1n_1  + e_{13}n_3n_1\right)
\notag\\ &\quad  
+\frac{\eps}{\Gamma}\frac{\xi q}{2}n_1n_1(q_t+v_kq_{,k})
+\frac{\eps}\Gamma [\xi\lambda_1 (1-q^2)-\zeta]\Delta\chi q (n_1 n_1-1/2),\lb{be:T_ten11}
\end{align}
\begin{align}
\tilde{T}_{13}=&
\alpha_1 \left(
\eps n_1n_1 e_{11}n_1n_3+
n_1n_3 e_{13}n_1n_3+
n_3n_1 e_{31}n_1n_3+
\eps n_3n_3 e_{33}n_1n_3
\right)
\notag \\ & \quad 
+ \alpha_2 N_1n_3 + \alpha_3N_3n_1 
+\alpha_4 e_{13} 
\notag \\ &\quad
+\alpha_5 \left(\eps e_{11}n_1n_3+ e_{13}n_3n_3\right) 
+\alpha_6 \left(e_{31}n_1n_1+\eps e_{33}n_3n_1\right)
\notag\\ & \quad 
+\frac{\eps}{\Gamma}\frac{\xi q}{2}n_1n_3(q_t+v_kq_{,k})
+ \frac{\eps}\Gamma[\xi\lambda_1 (1-q^2)-\zeta]\Delta\chi q n_1 n_3,\lb{be:T_ten13}
\end{align}
\begin{align}
\tilde{T}_{31}=&
\alpha_1 \left(
\eps n_1n_1 e_{11} n_3n_1 +
n_1n_3 e_{13} n_3n_1 +
n_3n_1 e_{31} n_3n_1 +
\eps n_3n_3 e_{33} n_3n_1 
\right)
\notag\\&\quad
+ \alpha_2 N_3n_1 + \alpha_3N_1n_3 
+\alpha_4 e_{31} 
\notag\\ & \quad 
+\alpha_5 \left(e_{31}n_1n_1+\eps e_{33}n_3n_1\right) 
+ \alpha_6 \left(\eps e_{11}n_1n_3+e_{13}n_3n_3\right)
\notag\\ & \quad 
+\frac{\eps}{\Gamma}\frac{\xi q}{2}n_3n_1(q_t+v_kq_{,k})
+   \frac{\eps}\Gamma[\xi\lambda_1 (1-q^2)-\zeta]\Delta\chi q n_3 n_1,\lb{be:T_ten31}
\end{align}
\begin{align}
\tilde{T}_{33}=&
\alpha_1\left( 
\eps n_1n_1 e_{11}n_3n_3 + 
n_1n_3 e_{13}n_3n_3 +
n_3n_1 e_{31}n_3n_3 +
\eps n_3n_3 e_{33}n_3n_3
\right)
\notag\\ & \quad 
+ \alpha_2 N_3n_3 + \alpha_3N_3n_3 
+\alpha_4 \eps e_{33} 
\notag\\ & \quad 
+\alpha_5 \left(e_{31}n_1n_3 + \eps e_{33}n_3n_3\right)
+\alpha_6 \left(e_{31}n_1n_3 + \eps e_{33}n_3n_3\right)
\notag\\ & \quad 
+\frac{\eps}{\Gamma}\frac{\xi q}{2}n_3n_3(q_t+v_kq_{,k})
+   \frac{\eps}\Gamma[\xi\lambda_1 (1-q^2)-\zeta]\Delta\chi q (n_3 n_3-1/2), \lb{be:T_ten33}
\end{align}
\end{subequations}
\begin{subequations}
\begin{align}
e_{11}=&\partial_1 v_1,&
\omega_{11}=&0,\\
e_{13}=&\frac12\left(\partial_3 v_1 + \eps^2\partial_1 v_3\right),&
\omega_{13}=&\frac12\left(\partial_3 v_1 - \eps^2\partial_1 v_3\right),\\
e_{31}=&\frac12\left(\eps^2\partial_1 v_3 + \partial_3 v_1\right),&
\omega_{31}=&\frac12\left(\eps^2\partial_1 v_3 - \partial_3 v_1\right),\\
e_{33}=&\partial_3 v_3 ,&
\omega_{33}=&0,
\end{align}
\end{subequations}

\begin{subequations}
\begin{align}
N_1=& \varepsilon\pt n_1 + \varepsilon v_j\pj n_1-\frac12\partial_3 v_1 n_3
+\varepsilon^2 \frac12 \partial_1 v_3 n_3,\\
N_3=& \varepsilon\pt n_3 + \varepsilon v_j\pj n_3-\varepsilon^2\frac12\partial_1 v_3 n_1
+\frac12\partial_3 v_1 n_1.
\end{align}
\end{subequations}

In turn, at the substrate $x_3=0$, the non-dimensional boundary conditions are
\begin{subequations}
\begin{align}
v_1&=0, \quad v_3=0,\\
n_3 &=\cos\theta_1,\\
q &= q_1,
\end{align}
\end{subequations}
and at the free surface, $x_3=\eta(x_1,t)$, they are
\begin{subequations}
\begin{align}
\eta_t &= v_3-v_1 \partial_1\eta,\\
-p+
\frac{\eps}{(1+\eps^2\eta_{1,1}^2)}
\left[(\eps^2T_{11}^E+\tilde T_{11})\eta_{1,1}^2 
\right.\qquad &
\notag \\
\left.
-(\eps T_{13}^E+\tilde T_{13})\eta_{1,1}
-(\eps T_{31}^E+\tilde T_{31})\eta_{1,1}
\right.\qquad &
\notag \\
\left.
+(T_{33}^E+\tilde T_{33})\right]
&= 
\frac{\eta_{,11}}{(1+\eps^2\eta_{1,1}^2)^{3/2}} ,
\\
-\eps\eta_{1,1}(\eps^2 T_{11}^E+\tilde T_{11}) - \eps^2\eta_{1,1}^2 (\eps T_{13}^E+\tilde T_{13}) 
\qquad &
\notag \\
+ (\eps T_{31}^E+\tilde T_{31}) + \eps\eta_{1,1}(T_{33}^E+\tilde T_{33}) &=0,\\
\frac{-\eps \eta_{1,1} n_1+n_3}{(1+\eps^2\eta_{1,1}^2)^{1/2}}&=\cos\theta_2,\\
q&=q_2,
\end{align}
\lb{NSCond}
\end{subequations}
where in the normal stress equation we have used again our earlier choice for $\eps$
in \eqref{new:eps}. 

\section{Derivation of the thin-film model for the active Eriksen-Leslie-Parodi theory}

\subsection{Governing equations}
In this appendix, we give a brief account of the derivation of the thin-film model
for the Eriksen-Leslie-Parodi theory augmented by activity terms. Conventions and notations carry over from the main text. The conservation of mass, 
linear and angular momentum balance equations are given by 
\begin{align}
0 &=\pai v_i, \label{cont}\\
0 &=-\pai p -\pj\left(\ppjnk W\,\pai n_k\right) + \pj\tilde T_{ij}
		\label{mom} \\
0 &= h_i - \gamma_1 N_i - \gamma_2 e_{ij} n_j+\lambda^{ELP}_1\Delta\chi^{ELP} n_i,
\label{mol}		
\end{align}

The bulk free energy density $W$ is
\begin{align}
2W=&K_1\left(\nabla\cdot\bn\right)^2+K_2\left(\bn\cdot\mbox{curl}\,\bn\right)^2
\notag\\
&\qquad +K_3\left(\bn\times\mbox{curl}\,\bn\right)^2+(K_2+K_4)[\tr(\nabla
  n)^2-(\nabla\cdot\bn)^2]
\label{W}
\end{align}
The parameters $K_1$, $K_2$ and $K_3$ are the splay, twist and bend elastic 
moduli (see de Gennes \& Prost \cite{GP93}), and $K_2+K_4$ is the saddle-splay constant.
Notice that in the case of strong anchoring, the final term does not contribute
to the governing equations \cite{Lin2013}, and that in 2D, there is also no twist term. 
Again, we will assume that all the $K_1=K_2=K_3\equiv K$ and $K_4=0$.
This assumption is discussed for liquid crystals in section 3.1.3.2 of \cite{GP93}),
and we use it here for simplification; see also \cite{VJP05}.  Notice that
under this assumption the elastic energy is reduced to the Dirichlet energy
(see also equation (4) in \cite{Lin2013})
\begin{equation}\label{WasDirEn}
2W=K \pk n_{i} \pk n_{i}. 
\end{equation}
In turn, the rate of change of the director with respect to the background
fluid $N_i$ is defined as in \rf{Ni}.
The molecular field is given by 
\be
h_i=\gamma n_i -\frac{\delta W}{\delta n_i}\label{molfield}
\ee
where $\gamma$
appears as a Lagrange multiplier in the variational formulation
to satisfy the condition $n_in_i=1$ and may in general
depend on $x_i$ and $t$.  

Here and in the following sections, we consider 
several stress tensors. The total stress tensor 
is given by
\begin{equation}\label{totstress}
T_{ij}=-p\delta_{ij}+ T^E_{ij}  +   \tilde T_{ij},
\end{equation}
where the Eriksen-Leslie tensor is
\begin{equation}\label{teij}
T^E_{ij} = -\ppink W\,\pj n_k
\end{equation}
and the extra stress tensor is
\begin{align}
\tilde T_{ij} &= \alpha_1 n_kn_p e_{kp}n_in_j + \alpha_2 N_in_j + \alpha_3N_jn_i 
+\alpha_4 e_{ij} \notag\\
& \quad + \alpha_5 e_{ik}n_kn_j + \alpha_6 e_{jk}n_kn_i
+\zeta^{ELP}\Delta\chi^{ELP} n_in_j.
\label{extrastressl}
\end{align}
We remark that in some of the literature, e.g. \cite{rey00}, $T_{ij}$ includes an
additional term $-W\delta_{ij}$, which, however, amounts to a redefinition
of the pressure \cite{Lin2013}. As in the nematic system \rf{cont1}-\rf{mol1} we introduced
in system \rf{cont}-\rf{mol} two active parameters $\lambda^{ELP}_1$ and $\zeta^{ELP}$.

\paragraph{Director field boundary conditions.}

At the substrate $x_3=0$, the strong anchoring condition reads
\be
n=\sin \theta_1 \, e_1 + \cos \theta_1 \,e_3, \quad \label{bc1-dir}
\ee
where $e_1$ and $e_3$ are the canonical unit vectors,
and at the free surface $x_3=\eta(x_,t)$, we have analogously
\be
n = \cos\theta_2 \,\nu + \sin\theta_2 \,t. \label{bc2a-dir}
\ee

\paragraph{Flow field and stress boundary conditions.}

For the boundary conditions of the flow field we 
assume at the substrate $x_3=0$ no-slip and impermeability, respectively 
\be
v_1=0,\qquad v_3=0, \label{bc1flow}
\ee
and at the interface $x_3=\eta(x_1,t)$, we have the kinematic condition
\begin{equation}\label{eqn:kincond}
\partial_t \eta=v_3- v_1 \partial_1\eta.
\end{equation}
The interfacial stress boundary condition is 
\be
\nu_i T_{ij}=-g_0 \pai\nu_i \nu_j, \label{bc2flow}
\ee
or in components
\bea
\nu_i T_{ij}\nu_j &=& -g_0 \pai\nu_i, \label{bcn-flow}\\
\nu_i T_{ij} t_j &=& 0.		\label{bct-flow}
\eea
From \eqref{totstress}, we obtain for the normal 
bulk stress
\begin{eqnarray}
\nu_i T_{ij}\nu_j&=&-p +\nu_i  T^E_{ij}\nu_j+\nu_i \tilde T_{ij}\nu_j,
\label{tijvivj}\\
\nu_i  T^E_{ij}\nu_j
&=&
-\nu_i\ppink W\,\pj n_k\nu_j,
\\
\nu_i \tilde T_{ij}\nu_j&=& 
\alpha_1 n_kn_p e_{kp}\nu_in_in_j\nu_j + \alpha_2 \nu_iN_in_j\nu_j + \alpha_3\nu_i N_jn_i\nu_j 
\nn\\&& 
+\alpha_4\nu_i e_{ij}\nu_j 
+ \alpha_5\nu_i e_{ik}n_kn_j\nu_j + \alpha_6\nu_i e_{jk}n_kn_i\nu_j
\nn\\&&
+\zeta^{ELP}\Delta\chi^{ELP}\nu_in_in_j\nu_j.
\end{eqnarray}
Similarly, we obtain for the tangential boundary condition 
\begin{eqnarray}
\nu_i T_{ij}t_j&=&\nu_i  T^E_{ij}t_j+\nu_i \tilde T_{ij}t_j,\\
\nu_i  T^E_{ij}t_j
&=&
-\nu_i\ppink W\,\pj n_kt_j,
\\
\nu_i \tilde T_{ij}t_j&=& 
\alpha_1 n_kn_p e_{kp}\nu_in_in_jt_j + \alpha_2 \nu_iN_in_jt_j + \alpha_3\nu_i N_jn_it_j 
\nn\\&& 
+\alpha_4\nu_i e_{ij}t_j 
+ \alpha_5\nu_i e_{ik}n_kn_jt_j + \alpha_6\nu_i e_{jk}n_kn_it_j 
\nn \\&&
+\zeta^{ELP}\Delta\chi^{ELP}\nu_in_in_j t_j.
\end{eqnarray}


\subsection{Thin-film approximation}

Using the same nondimensionalisation \rf{sc_b}--\rf{new:eps} as before together with $W=\wscale \bar W$ and
\be
\wscale=\frac{K }{\eps^2 L^2}\label{new:Gamma},
\ee
we obtain that the  non-dimensional bulk free energy becomes
\be
2\bar W=(\pthr n_3)^2 + (\pthr n_1)^2 + O(\eps^2).
\label{sW0}
\ee
The scaled  molecular field is then to leading order
\bea
\bar h_1&=& n_1 + \pthr^2 n_1 + O(\eps^2)\label{molfield0a}\\
\bar h_3&=& n_3 + \pthr^2 n_3 + O(\eps^2)\label{molfield0b}
\eea

Further on, we introduce the dimensionless parameters
\be
\bar\alpha_i=\alpha_i/\mu,\ \bar\gamma_i=\gamma_i/\mu,
\Delta\bar\chi^{ELP}=\frac{\eps L}{\mu U}\Delta\chi^{ELP},\ \bar\zeta^{EPL}=\zeta^{EPL}/\mu,
\ee
where $\mu$ is the kinematic viscosity. Then we nondimensionalize \rf{mol} to obtain (upon neglecting lower order terms) that
\bea
0 &=&\left(\frac{K}{\eps^2 L^2}\right) \,\bar h_1 - 
\left(\frac{\alpha_2}{\mu}\right)\left(\frac{\mu U}{\eps L}\right)\, n_3\pthr
\bar v_1+\frac{\mu U}{\eps L}\bar\lambda^{ELP}_1\Delta\bar\chi^{ELP} n_1
\label{smola}\\
0 &=&\left(\frac{K}{\eps^2 L^2}\right)\, \bar h_3 - 
\left(\frac{\alpha_3}{\mu}\right)\left(\frac{\mu U}{\eps L}\right)\, n_1\pthr \bar v_1+\frac{\mu U}{\eps L}\bar\lambda^{ELP}_1\Delta\bar\chi^{ELP} n_3
\label{smolb}
\eea
If $K/\eps^2 L^2 \gg \mu U/\eps L$ the flow field decouples from  the director field in these equations. Therefore, we require  the case of weak elasticity $K/\eps\mu U L=O(1)$, so that all three terms in each of the equations \eqref{smola} and \eqref{smolb} remain.

The scale $P$ for the pressure is obtained, as before, by balancing 
it with the dominant viscous contributions in the horizontal momentum 
equation \rf{mom} (i.e. for $i=1$).  
We drop the overbars from this point onwards and introduce
$\theta$ as in \eqref{n_Rep}. The leading order bulk equations then are
\begin{subequations}
\begin{align}
0&=v_{1,1}+v_{3,3}, \label{xcont}\\
0&=-p_{,1}+\frac{1}{2}(v_{1,3}f_A(\theta))_{,3}+\zeta^{ELP}\Delta\chi^{ELP}
(\sin(2\theta))_{,3},
\\
0&=-p_{,3},
\\
0  &= \gamma \sin\theta + K(\sin\theta)_{,33}+\frac{1}{2}(\gamma_1-\gamma_2) 
v_{1,3}\cos\theta+\lambda^{ELP}_1\Delta\chi^{ELP}\sin\theta,\\
0  &= \gamma \cos\theta + K(\cos\theta)_{,33}-\frac{1}{2}(\gamma_1+\gamma_2) 
v_{1,3}\sin\theta+\lambda^{ELP}_1\Delta\chi^{ELP}\cos\theta,
\end{align}
with the Lagrange parameter $\gamma$ and
\begin{align}
f_A(\theta)&=(\alpha_1/2)\sin^2(2\theta)+(\alpha_5-\alpha_2)\cos^2\theta+(\alpha_3+\alpha_6)\sin^2\theta+\alpha_4.
\end{align}
\lb{EM_pol}
\end{subequations}
Notice that in the above all terms in 
$\pj\left(\ppjnk W\,\pone n_k\right)$ 
are of order $\eps$ or smaller and hence do not contribute.
The leading order boundary conditions are as follows:  At $x_3=0$, we have 
\begin{subequations}
\begin{align}
v_1&=0, \quad v_3=0,   \label{xbc1flow}\\
\theta&=\theta_1,
\end{align}
\end{subequations}
and at the free surface, $x_3=\eta(x_1,t)$, 
\begin{subequations}
\begin{align}
\eta_t &= v_3-v_1 \partial_1\eta, \label{eqn:xkincond}\\ 
-p&= \eta_{,11} ,\\
\frac{1}{2}v_{1,3}f_B(\theta_2) &=-\frac{\zeta^{ELP}\Delta\chi^{ELP}}{2}\sin(2\theta_2)
,\\
\theta&=\theta_2,
\end{align}
\end{subequations}
where we have defined
\begin{align}
f_B(\theta)&=(\alpha_1/2)\sin^2(2\theta)
+(\alpha_6-\alpha_3)\cos^2\theta+(\alpha_2+\alpha_5)\sin^2\theta+\alpha_4.
\lb{EM_pol1}
\end{align}
Similarly as it was done in~\cite{Lin2013} for its passive counterpart, the system \rf{EM_pol}--\rf{EM_pol1} can be partly integrated to yield the active thin-film model \rf{LEP_LM} based on the Leslie-Erickson-Parodi theory. In particular, the mass conservation relation \rf{kin} can be derived from \eqref{xcont}, \eqref{eqn:xkincond} and \eqref{xbc1flow}.

\end{appendix}

\bibliographystyle{abbrvnat}
\bibliography{articles}

\end{document}

%% file: macros-cell.tex
\def\tr{{\rm tr}}
\newcommand{\lb}{\label}
\newcommand{\mR}{\mathbb{R}\,}
\newcommand{\mRtwo}{\mathbb{R}^2}
\newcommand{\mZ}{\mathbb{Z}\,}
\newcommand{\id}{I}

\newcommand{\rf}[1]{(\ref{#1})}

\newcommand{\lubp}{\varepsilon}
\newcommand{\bl}{\begin{linenomath}}
\newcommand{\el}{\end{linenomath}}

\newcommand{\bn}{{\bf n}}

\newcommand{\eps}{\varepsilon}

\newcommand{\const}{\mathrm{const}}

\newcommand{\pai}{\partial_i}

\newcommand{\pj}{\partial_j}
\newcommand{\pk}{\partial_k}

\newcommand{\ppink}{\partial_{\partial_in_k}}
\newcommand{\ppjnk}{\partial_{\partial_jn_k}}

\newcommand{\pt}{\partial_t}

\newcommand{\pone}{\partial_1}
\newcommand{\pthr}{\partial_3}

\newcommand{\be}{\begin{equation}}
\newcommand{\ee}{\end{equation}}
\newcommand{\bea}{\begin{eqnarray}}
\newcommand{\eea}{\end{eqnarray}}
\newcommand{\sbea}{\begin{subequations}\begin{eqnarray}}
\newcommand{\seea}{\end{eqnarray}\end{subequations}}
\newcommand{\beas}{\begin{eqnarray*}}
\newcommand{\eeas}{\end{eqnarray*}}
\newcommand{\nn}{\nonumber}
\newcommand{\ees}{\end{equation*}}
\newcommand{\bes}{\begin{equation*}}
